\documentclass[12pt]{article}
\usepackage{amssymb}
\usepackage{amsmath}
\usepackage{amsthm}
\usepackage{amsfonts}
\usepackage[dvips]{graphicx}
\usepackage{color}

\usepackage{bbm,verbatim}

\alph{footnote}
\title {The expected area of the filled planar\\ Brownian loop is $\frac{\pi}{5}$}
\author {Christophe Garban \footnote {Department of Mathematics, Cornell University and ENS Paris.
Email: christophe.garban@ens.fr}
 \and  Jos\'e A. Trujillo Ferreras \footnote {Department of Mathematics, Cornell University. Email:
jatf@math.cornell.edu}}

\newtheorem{theorem}{Theorem}[section]
\newtheorem{lemma}[theorem]{Lemma}

\theoremstyle{definition}

\theoremstyle{remark}

\numberwithin{equation}{section}

\newcommand{\Z}{{\mathbb Z}}

\newcommand{\R}{{\mathbb{R}}}
\newcommand{\C}{{\mathbb C}}
\renewcommand{\H}{{\mathbb H}}

\newcommand{\slemeas}{\mu ^{{\rm sle}}}
\newcommand{\gs}{\gamma ^*}
\newcommand{\eps}{\varepsilon}

\newcommand{\f}{\frac}

\newcommand{\Prob} {{\bf P}}
\newcommand{\E}{{\bf E}}
\renewcommand{\aa}{\alpha}

\def \me {\text{m}_{\varepsilon}}
\def \mei {\text{m}_{\varepsilon}^{-1}}

\def \A {\mathcal{A}}

\def \mubb {\mu^{\text{bub}}}
\def \musle {\mu^{\text{sle}}}
\def \ad {\musle(A|\gs=1)}
\def \Dp {\mathbb{D}^+}
\def \z {{z_0}}

\begin{document}
\maketitle
\alph{footnote}
\date{}
\begin{abstract}
Let $B_t, 0 \le t \le 1$ be a planar Brownian loop (a Brownian
motion conditioned so that $B_0=B_1$). We consider the compact
hull obtained by filling in all the holes, i.e. the complement of
the unique unbounded component of $\C\setminus B[0,1]$. We show
that the expected area of this hull is $\pi / 5.$ The proof uses,
perhaps not surprisingly, the Schramm Loewner Evolution ($SLE$).
Also, using Yor's result \cite{yor} about the law of the index of
a Brownian loop, we show that the expected areas of the regions of
index (winding number) $n\in\Z\setminus\{0\}$ are $\frac{1}{2\pi
n^2}\,$. As a consequence, we find that the 
expected area of the region of 
index zero inside the loop is $\pi / 30$; this 
value could not be obtained directly using 
Yor's index description. 
\end{abstract}

\begin{figure}
\begin{center}
\includegraphics[width=0.65\textwidth]{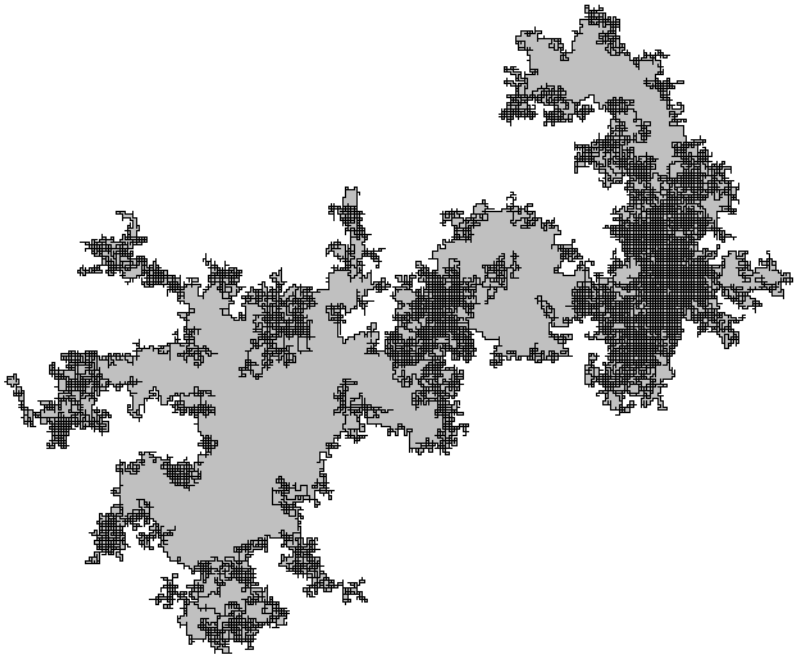}
\end{center}
\caption{Random walk loop of 50000 steps and corresponding hull.}
\end{figure}

\section{Introduction}
In the abundant literature about planar Brownian motion, there are
certainly results dealing with the question of area.
Paul L\'evy's stochastic area formula  describing the algebraic area
``swept'' by a Brownian
motion will likely come to the mind of many readers. Our result, however,
 is very different from this classical theorem,
firstly because Levy's area is a signed area, but mainly because
of the following : in order to apprehend L\'evy's area it
is enough to follow the Brownian curve locally without paying
attention to the rest of the curve. In our case, one
needs to consider the curve globally.

Aside from the fact that the question we address is a very natural
one for Brownian motion, we
have been motivated by related results in the Physics literature.
In \cite{cardy}, using methods of conformal field theory, Cardy
has shown that the ratio of the expected area enclosed by a
self-avoiding polygon of perimeter $2n$ to the expected squared radius of
gyration for a polygon of perimeter $2n$ converges as $n$ goes to
infinity to $4 \pi / 5$. We note that self-avoiding
polygons are supposed to have the same asymptotic shape
as filled Brownian loops (see, for example, \cite{richard} and
references therein).  However, studying this relationship is hard
basically for the following reason. The boundary of the Brownian loop is
of $SLE _{8/3}$-type, but, unfortunately, there does not exist a good
way of ``talking about the length'' of $SLE$ curves at this moment.

Our result gives interesting information regarding the
Brownian loop soups introduced in \cite{BLS}.
This conformally invariant object plays
an important role in the understanding and description of $SLE$
 curves (see, e.g. \cite{BLS,RQ,CR}).
It can be viewed as a Poissonian cloud (of intensity $c$)
of {\em filled} Brownian loops in
subdomains of the plane. Among other things, it is announced in \cite{W}
that the
dimension of the set of points in the complement of the loop soup (i.e. the
points that are in the inside of no loop) can be shown to be equal to
$2-c/5$, using consequences of the restriction property.
A detailed proof of this
statement has never been published, and in fact, our result implies the
corresponding first moment estimate (i.e. the mean number of balls of radius
$\eps$ needed to cover the set). The other arguments needed to derive the
result announced in \cite{W} will be detailed in \cite{T}.

Let us make precise what we mean by area enclosed by a Brownian loop.
Let $B$ denote a Brownian bridge in $\C$ of time duration 1.
I.e. the law of $B_t,\, 0 \le t \le 1$ is the same as the law of
$W_t - t W_1, \, 0 \le t \le 1$, where $W$ is just a standard
Brownian motion in $\C$.
 $\C\setminus
B[0,1]$, i.e. the complement of the path, has a unique infinite
connected component $H$. The hull $T$ generated by the Brownian
loop is by definition $\C\setminus H$. Let $\cal{A}$ be the random
variable whose value is the area of $T$. In this paper, we will
prove

\begin{theorem}\label{main}
\[
\E(\A)=\frac{\pi}{5}.
\]
\end{theorem}

We would like to explain now how this result is related 
to the problem of windings of a Brownian loop. 
In \cite{yor}, Yor gave an explicit
formula for the law of the index of a Brownian loop around a fixed point $z$.
A point with a non-zero index has to be inside the loop. Using
this fact, it is almost possible to describe the probability that a
point is inside the loop, modulo the problem of the
zero index; indeed, there are some regions inside the Brownian loop
which are of index zero. It seems hard to control the influence of
these zero-index points inside the curve. In the last section,
using our main result, theorem~\ref{main},
combined with the law of the index given by Yor,
we find that the expected area of the
set of points inside the loop that have index zero is $\frac{\pi}{30}$. 
We also compute the expected areas of the regions of index
$n\in\Z\setminus \{0\}$. 

In \cite{winding}, using physics methods, 
Contet, Desbois and Ouvry obtained the values of the 
expected areas for the non-zero regions. 
In their paper, they noted the different nature of the $n=0$ sector (the
points in the plane of zero index) and emphasized that ``it would
be interesting to distinguish in the $n=0$ sector, curves which do
not enclose the origin from curves which do enclose the origin but
an equal number of times clockwise and anticlockwise''. Their values
in the case $n\neq0$ agree with our results; 
they argue that the
0-case cannot be treated within the scope of their analysis.

From a probabilistic viewpoint, it also appears that 
usual techniques for Brownian motion are not
strong enough to obtain the expected area of 
the Brownian loop or the expected area of the 0-index 
region inside the Brownian loop. However, the computation of 
the expected area of the $n$-index region for $n \ne 0$ 
was within reach using the result of Yor. To our knowledge
this computation had not been carried out in a mathematical way
before. 

Let us briefly explain why usual 
techniques for Brownian motion 
seem 
unable to tackle the problem of the expected area of the 
Brownian loop. Basically, the enclosed area
depends only on the boundary of the hull generated by the Brownian
loop. The frontier of the Brownian loop concerns only a small
subset of the time duration $[0,1]$. In some sense, on certain
time-intervals, the enclosed area does not depend much on the
behavior of the Brownian motion. So, this problem needs a good
description of the frontier of a Brownian loop. Recently, Lawler,
Schramm and Werner proved a conjecture of Mandelbrot that the
Hausdorff dimension of the Brownian frontier is $4/3$. For this
purpose they used the value of intersection exponents computed
with the help of $SLE$ curves, see for instance \cite {boundary}
and references therein. The description of the Brownian frontier
via $SLE$ can be done in a slightly different way using the
conformal-restriction point of view, see \cite {CR}. We will use
this approach, and so will present to the reader the facts needed
about conformal restriction measures in the next section.

 The paper gives another striking example of a simple result
concerning planar Brownian motion that seemed out of reach using
the usual stochastic calculus approach, but that can be derived
using conformal invariance and $SLE$. For a thorough account on
$SLE$ processes, see \cite{book,stflour}.

\section{Preliminaries}\label{preliminaries}
Conformal restriction measures in $\H$ are measures supported on
the set of closed subsets $K$ of $\H$ such that
$\overline{K}\cap\R=\{0\}$, $K$ is unbounded and $\H\setminus K$
has two infinite connected components, that satisfy the conformal
restriction property : for all simply connected domains
$H\subset\H$ such that $\H\setminus H$ is bounded and bounded away
from the origin, the law of $K$ conditioned on $K\subset H$ is the
law of $\Phi(K)$, where $\Phi$ is any conformal transformation
from $H$ to $\H$ preserving $0$ and $\infty$ (this law doesn't
depend of the choice of $\Phi$). It is proved in \cite {CR} that
there is only one real parameter family of such restriction
measures, $\Prob_{\aa}$ where $\aa\geq 5/8$. These measures are
uniquely described by the following property : for all closed
$A$ in $\H$ bounded and bounded away from
$0$,
\begin{equation}\label{restriction}
\Prob_{\aa}[K\cap
A=\emptyset]=\Phi_A^{\,'}(0)^{\aa}\,,
\end{equation}
where $\Phi_A$ is a conformal transformation from $\H\setminus A$
into $\H$ such that $\Phi_A(z)/z\rightarrow 1$, when
$z\rightarrow\infty$. To aid with the notation for the rest of the
paper whenever we write $\Phi _A$ we will be assuming that we have
chosen the translate with the additional property $\Phi _A (0) =
0$. $\Prob_{5/8}$ is the law of chordal $SLE_{8/3}$, and $\Prob_1$
can be constructed by filling the closed loops of a Brownian
excursion in $\H$ (Brownian motion started at $0$ conditioned to
stay in $\H$). An important property of these conformal
restriction measures is that using two independent restriction
measures $\Prob_{\aa_1}$ and $\Prob_{\aa_2}$, we can construct
$\Prob_{\aa_{1}+\aa_2}$ by filling the ``inside'' of the union of
$K_1$ and $K_2$. This ``additivity'' property and the construction
of $\Prob_{5/8}$ and $\Prob_1$ give the good description of the
Brownian motion in terms of $SLE$ curves, namely, 8 $SLE_{8/3}$
give the same hull as 5 Brownian excursions.

Since, we want to describe the
boundary of loops of time duration 1, we will first create
loops with the use of the infinite hulls described above.
Restriction measures are conformally invariant (Brownian
excursion, $SLE_{8/3}$,..), so we had better use conformal maps.
There is obviously no conformal map which sends both $\infty$ and 0
to 0, so the natural idea is to consider a M\"{o}bius
transformation preserving $\H$ which maps 0 to 0, and $\infty$ to
$\eps$. We can choose
\[
\me(z)=\frac{\eps z}{z+1}
\]
\[
\mei(z)=\frac{z}{\eps-z}\,.
\]
The limit when $\eps$ goes to zero of the measures $\me(\Prob_1)$ is
the dirac measure at $\{0\}$. The good renormalization to keep
something interesting is in $\eps^2$. Hence, we define the Brownian
bubble measure in $\H$ as :
\[
\mubb=\lim_{\eps\rightarrow 0}\frac{1}{\eps^2}\me(\Prob_1)\,.
\]

This measure was introduced in \cite{CR}, and it is an
important tool for studying the link between $SLE$ curves and the Brownian
loop soup (see \cite {BLS}). It was already noted in \cite{CR,saw}, as an
easy consequence of the ``additivity'' property described above,
that
\[
\frac{5}{8}\mubb=
\frac{5}{8}\lim_{\eps\rightarrow 0}\frac{1}{\eps^2}\me(\Prob_1)=
\lim_{\eps\rightarrow 0}\frac{1}{\eps^2}\me(\Prob_{5/8})\,.
\]
The last measure can be seen as an infinite measure on
``$SLE_{8/3}$ loops", let us call this measure $\musle$. Recall, that
we are interested in a Brownian loop of time duration 1.
We have the following time
decomposition for $\mubb$, (see \cite{BLS},\cite{book})
\begin{equation}
\label{decomp}
 \mubb=\int_0^{\infty}\frac{dt}{2
\,t^2}\Prob_t^{\text{br}}\times\, \Prob_t^{\text{exc}}\,,
\end{equation}
 where  $\Prob_t^{\text{br}}$ is the law of
a {\em one-dimensional}
Brownian bridge of time duration $t$, and
$\Prob_t^{\text{exc}}$ is the law of an It\^{o} Brownian excursion
re-normalized to have time $t$.
$\Prob_t^{\text{br}}\times\, \Prob_t^{\text{exc}}$ is the law of an
$\H$-Brownian bridge of time duration $t$, by considering
the one dimensional bridge as the
$x$ coordinate of the curve, and the excursion as the
$y$ coordinate. Unfortunately, it is hard to compute fixed-time
quantities with $SLE$ techniques. Thus, we will compute a ``geometric quantity"
using $SLE_{8/3}$, and then extract $\E(\A)$ from this
geometric value by using the relation $\mubb=8/5\musle$ and
the decomposition ~\ref{decomp}.

Let us explain in a few words why we need to deal with Brownian
bridges in $\H$ and cannot work directly with bridges in $\C$.
The
underlying idea is the fact that one needs to choose a starting point
on the boundary of the Brownian loop for the $SLE$ loop
representation. A natural choice is the (almost surely) unique
lower point, this is why we are interested in $\H$ quantities.
So let $\A^\H$ be the random variable giving the
area of an $\H-$Brownian bridge
of time duration one. Working with $\A^\H$ will turn out not to be
a problem since, as the reader might already suspect, the random variables
$\A$ and $\A ^\H$ have the same law.

For the geometric quantity, we could choose to compute $\int A(\gamma)
d\musle$, where $A(\gamma)$ is the area enclosed in $\H$ by the
``curve" $\gamma$, but this integral is infinite.
Let $\gs$ be the radius of the curve $\gamma$, that is, $\gs=\sup_{0\leq t \leq
t_\gamma}|\gamma(t)|$. We may consider the ``expected" area under the law
$\musle$ ``conditioned" on $\gs =1$. Here, $\musle$ is not a
probability measure so the term ``expected value" is not correct,
and the conditioning is on a set of $\musle-$measure equal to 0.
But we
have the following rigorous definition :
\begin{equation}\label{Dpsle}
\ad=\lim_{\delta\downarrow 0}\frac{\int A(\gamma)1_{\{\gs\in
[1,1+\delta)\}} d\musle}{\musle\{\gs\in[1,1+\delta)\}}\,.
\end{equation}
Using $\musle=5/8\mubb$, we can write in the same way
:
\begin{equation}\label{Dpbb}\ad=\lim_{\delta\downarrow
0}\frac{\int A(\gamma)1_{\{\gs\in [1,1+\delta)\}}
d\mubb}{\mubb\{\gs\in[1,1+\delta)\}}\,.
\end{equation}
Thus, $\ad$
represents at the same time the ``expected" area of an $SLE_{8/3}$
loop conditioned to touch the half circle of radius one and the
expected area of a Brownian bubble with the same conditioning.
With the use of the restriction property for $SLE_{8/3}$, we will
be able to compute in the last section $\ad$. Before, in the coming
section, we will find the relationship between
$\E(\A)$ and
$\ad$.

\section{Extraction of $\E(\A)$ from $\ad$}

In this section we will prove the following
\begin{lemma}\label{extraction}
\[
\E ({\cal A}) = 2 \ad.
\]
\end{lemma}

\begin{proof}
First of all, by using the definition of $\mubb$
in terms of $\lim_{\eps\downarrow 0} \frac{1}{\eps^2}\me(\Prob_1)$
and the restriction property of $\Prob_1$, it is easy to show that
$\mubb\{\gs\geq r\}=\frac{1}{r^2}$, hence
\[
\mubb\{\gs\in[1,1+\delta)\}=1-1/(1+\delta)^2=2\delta+O(\delta^2),
\]
and thus, from \eqref{Dpbb}, we have
\begin{eqnarray*}
\ad & = & \lim_{\delta\downarrow 0} \frac{\int
A(\gamma)1_{\{\gs\in [1,1+\delta)\}}d\mubb}{\mubb\{\gs\in[1,1+\delta)\}} \\
& = & \lim_{\delta\downarrow 0}
\frac{\int_0^{\infty}\frac{dt}{2t^2}\E_t(A(\gamma)1_{\{\sup\limits_{0\leq
u\leq t}|\gamma(u)|\in[1,1+\delta)\}})}{2\delta+O(\delta^2)}.
\end{eqnarray*}
Here $\E_t$ is the expectation according to the law of an
$\H$-Brownian bridge in time $t$. By Brownian scaling
we have \[
\E_t(A(\gamma)\,1_{\{\sup\limits_{0\leq u\leq
t}|\gamma(u)|\in[1,1+\delta)\}})=t*\E_1(A(\gamma)\,1_{\{\sup\limits_{0\leq
u\leq
1}|\gamma(u)|\in[\frac{1}{\sqrt{t}},\frac{1}{\sqrt{t}}+\frac{\delta}{\sqrt{t}})\}}).
\]
Therefore :
\begin{eqnarray*}
\ad &=& \lim_{\delta\downarrow 0} \int_0^{\infty}
\frac{dt}{4\,t(\delta+O(\delta^2))}
\E_1(A\,1_{\{\gs\in[\frac{1}{\sqrt{t}},\frac{1}{\sqrt{t}}+\frac{\delta}{\sqrt{t}})\}})
\\ &=& \lim_{\delta\downarrow 0}
\int_0^{\infty}\frac{dt}{4\,t^{3/2}}
\frac{\Prob_1\{\gs\in[\frac{1}{\sqrt{t}},\frac{1}{\sqrt{t}}+\frac{\delta}{\sqrt{t}})\}}
{\frac{\delta}{\sqrt{t}}(1+O(\delta))}
\E_1(A\,|\gs\in[\frac{1}{\sqrt{t}},\frac{1}{\sqrt{t}}+\frac{\delta}
{\sqrt{t}})) \\ &=& \lim_{\delta\downarrow 0}
\frac{1}{2}\int_0^{\infty}du
(1+O(\delta))\frac{\Prob_1\{\gs\in[u,u+\delta u)\}}{\delta
u}\E_1(A\,|\gs\in[u,u+\delta u))\,,
\end{eqnarray*}
using the change of variables $u=\frac{1}{\sqrt{t}}\,$. Let $\eta_1$
be the density on $\R_+$ of the random variable $\gs$ under the
$\H$-Brownian bridge of time duration one. As for the one dimensional
bridge (law of the maximum of the bridge), this density decays
exponentially fast at infinity. Thus, we can interchange the limit and
the integral to obtain :
\begin{equation*}
\ad = \frac{1}{2}\int_0^{\infty}\eta_1(u)\E_1(A|\gs=u)du =
\frac{1}{2}\E(\A^\H)\,.
\end{equation*}
Hence, the proof of the lemma will be concluded as soon as we
establish
\[
\E(\A^\H)=\E(\A).
\]
 There is a (almost sure) one to one correspondence between $\C$-Brownian
bridges and $\H$-Brownian bridges. The idea is to
 start the Brownian loop from its lowest point. More
 precisely, if $B_t, 0\leq t\leq 1$ is a Brownian bridge in $\C$, with probability one,
 there is a unique $\bar{t}\in[0,1]$ such that $\text{Im}(B_{\bar{t}})\leq \text{Im}(B_t)$,
 for all $t\in [0,1]$. We associate to the Brownian Bridge $B_t$ the process
 $(Z_t)_{0\leq t \leq 1}$ in $\overline{\H}$, defined by this simple space-time translation :
 \begin{equation}
 \label{vv}
    Z_t = \left\{ \begin{array}{ll} B_{\bar{t}+t}-B_{\bar{t}} & ,0\leq t \leq 1-\bar{t} \,, \\
                       B_{\bar{t}+t-1}-B_{\bar{t}}& ,1-\bar{t} \leq t \leq 1\,. \end{array} \right.
                       \end{equation}
 Now, we have to identify the law of $Z_t$ with
 $\Prob_1^{\text{exc}}\times\Prob_1^{\text{br}}$. The real and imaginary
 parts of $B_t$ are two independent one-dimensional Brownian
 bridges. The law of the random variable $\bar{t}$ is independent of $\text{Re}(B_t)$,
 so in the space-time change \ref{vv}, $\text{Re}(Z_t)$ is still a
 one-dimensional bridge independent of the imaginary part of $Z_t$. $\text{Im}(Z_t)$
 has the law of a one-dimensional Brownian bridge viewed from its (almost sure) unique lowest point.
 By the Vervaat Theorem (see \cite{vervaat}), this gives the law
 of an It\^{o} excursion renormalized to have time one. Thus $Z_t$ has the
 law of an $\H$-Brownian bridge of time one. Our space-time transformation
 obviously preserves the area, hence $\E(\A^\H)=\E(\A)$.

\end{proof}

\section{Computation of $\ad$, and proof of theorem~\ref{main}}
In this section we prove lemma \ref{comp}, the proof provides a
good example of the use of standard techniques for $SLE_{8/3}$. We
have chosen to leave out some algebraic details in order to allow
the reader to focus on the main ideas.
\begin{lemma}\label{comp}
\[
\ad = \f \pi {10}.
\]
\end{lemma}

Note that theorem~\ref{main} follows immediately from this lemma
and lemma~\ref{extraction}.

\begin{proof}
Recall \eqref{Dpsle} :
\begin{equation}\label{ad}
\ad=\lim_{\delta\downarrow 0}\frac{\int A(\gamma)1_{\{\gs\in
[1,1+\delta)\}} d\musle}{\musle\{\gs\in[1,1+\delta)\}}\,.
\end{equation}
By using the definition $\musle=\lim_{\eps\downarrow
0}\frac{1}{\eps^2}\me(\Prob_{5/8})\,,$ we can rewrite \ref{ad} as :
\begin{equation}\label{eps}
\lim_{\delta\downarrow 0}\lim_{\eps\downarrow 0} \E_\eps
(A(\gamma)|\gs\in[1,1+\delta)), \end{equation} where $\E_\eps$ is
a more appealing notation for the expected value under the law of
$\me(\Prob_{5/8})$ (this law, in simpler words, is the law of a
chordal $SLE_{8/3}$ in $\H$ from 0 to $\eps$). Recall that
$A(\gamma)$ is the area of the bounded set in $\H$ enclosed by the
curve $\gamma$. $A(\gamma)$ can be written as $\int_\H 1_{\{z\text{
inside}\}} d{\rm A}(z)$, where \{$z$ inside\} means that $z$ is in the
component bounded by $\gamma$. Thus \eqref{eps} can be written as :
\begin{equation}
\label{eps2} \lim_{\delta\downarrow 0}\lim_{\eps\downarrow
0}\E_\eps\left(\int_{(1+\delta)\Dp}1_{\{z\text{
inside}\}}d{\rm A}(z)|\gs\in[1,1+\delta)\right)\,,
\end{equation}
where $\Dp$ is $\mathbb{D}\cap\H$. Since everything is nicely
bounded, we can interchange the limits and the integral. This
gives us :
\begin{equation}\label{goal}
\slemeas( A |\gs = 1) = \int_{\Dp} \lim _{\delta \downarrow 0}
\lim _{\eps \downarrow 0} \Prob_\eps \{ z \mbox{ inside } | \gs
\in [1, 1 + \delta) \} \,d{\rm A}(z).
\end{equation}

Therefore, what remains to be done is to compute, for a fixed z, the
``probability" that this $z$ is inside an ``$SLE_{8/3}$ loop"
conditioned to have radius exactly 1. So let us fix $z_0$ in
$\Dp$. Let $D_\eps$ ( resp $D_\eps^\delta$) denote the image under
$\text{m} ^{-1} _\eps (z) = z /(\eps - z)$ of the set $\{ z \in \H : |z|
\ge 1 \}$ (resp $\{ z \in \H : |z| \ge 1 +\delta\}$).

\begin{picture}(0,0)%
\includegraphics{area.pstex}%
\end{picture}%
\setlength{\unitlength}{3355sp}%
\begingroup\makeatletter\ifx\SetFigFont\undefined%
\gdef\SetFigFont#1#2#3#4#5{%
  \reset@font\fontsize{#1}{#2pt}%
  \fontfamily{#3}\fontseries{#4}\fontshape{#5}%
  \selectfont}%
\fi\endgroup%
\begin{picture}(8724,3446)(289,-6185)
\put(1525,-4442){\makebox(0,0)[lb]{\smash{\SetFigFont{10}{12.0}{\rmdefault}{\mddefault}{\updefault}{\color[rgb]{0,0,0}$z_0$}%
}}}
\put(2926,-4261){\makebox(0,0)[lb]{\smash{\SetFigFont{10}{12.0}{\rmdefault}{\mddefault}{\updefault}{\color[rgb]{0,0,0}$\mei$}%
}}}
\put(2326,-5011){\makebox(0,0)[lb]{\smash{\SetFigFont{10}{12.0}{\rmdefault}{\mddefault}{\updefault}{\color[rgb]{0,0,0}$1$}%
}}}
\put(2476,-5011){\makebox(0,0)[lb]{\smash{\SetFigFont{10}{12.0}{\rmdefault}{\mddefault}{\updefault}{\color[rgb]{0,0,0}$1+\delta$}%
}}}
\put(4351,-4411){\makebox(0,0)[lb]{\smash{\SetFigFont{10}{12.0}{\rmdefault}{\mddefault}{\updefault}{\color[rgb]{0,0,0}$\mei(z_0)$}%
}}}
\put(4801,-5086){\makebox(0,0)[lb]{\smash{\SetFigFont{10}{12.0}{\rmdefault}{\mddefault}{\updefault}{\color[rgb]{0,0,0}$0$}%
}}}
\put(3826,-5086){\makebox(0,0)[lb]{\smash{\SetFigFont{10}{12.0}{\rmdefault}{\mddefault}{\updefault}{\color[rgb]{0,0,0}$-1$}%
}}}
\put(5626,-3811){\makebox(0,0)[lb]{\smash{\SetFigFont{10}{12.0}{\rmdefault}{\mddefault}{\updefault}{\color[rgb]{0,0,0}$\Phi_{\eps}$}%
}}}
\put(7651,-4111){\makebox(0,0)[lb]{\smash{\SetFigFont{10}{12.0}{\rmdefault}{\mddefault}{\updefault}{\color[rgb]{0,0,0}$0$}%
}}}
\put(7276,-3511){\makebox(0,0)[lb]{\smash{\SetFigFont{10}{12.0}{\rmdefault}{\mddefault}{\updefault}{\color[rgb]{0,0,0}$\Phi_{\eps}(\mei(z_0))$}%
}}}
\put(7726,-6136){\makebox(0,0)[lb]{\smash{\SetFigFont{10}{12.0}{\rmdefault}{\mddefault}{\updefault}{\color[rgb]{0,0,0}$0$}%
}}}
\put(7426,-5386){\makebox(0,0)[lb]{\smash{\SetFigFont{10}{12.0}{\rmdefault}{\mddefault}{\updefault}{\color[rgb]{0,0,0}$\Phi_{\eps}^{\delta}(\mei(z_0))$}%
}}}
\put(1651,-5086){\makebox(0,0)[lb]{\smash{\SetFigFont{10}{12.0}{\rmdefault}{\mddefault}{\updefault}{\color[rgb]{0,0,0}$\eps$}%
}}}
\put(1501,-5086){\makebox(0,0)[lb]{\smash{\SetFigFont{10}{12.0}{\rmdefault}{\mddefault}{\updefault}{\color[rgb]{0,0,0}$0$}%
}}}
\put(6151,-5311){\makebox(0,0)[lb]{\smash{\SetFigFont{10}{12.0}{\rmdefault}{\mddefault}{\updefault}{\color[rgb]{0,0,0}$\Phi_{\eps}^{\delta}$}%
}}}
\end{picture}

We warn the reader that $\gamma$ will denote two different kinds
of curves in $\H$ : a curve from 0 to $\infty$, or a curve from
$0$ to $\eps$. Let $F_\eps$ be the event $\{ \gamma[0,\infty) \cap
D _\eps = \emptyset \}$, and, similarly, let $F _\eps ^\delta$ be
the analogous event for $D _\eps ^\delta$. Then,
\begin{equation*}
\Prob _\eps \{ z_0 \mbox{ inside } | \gs \in [1, 1 + \delta )\} =
\Prob_{5/8} \{\text{m} _\eps ^{-1} (z_0) \mbox{ is to the right of
$\gamma$ } | (F_ \eps)^c \cap F_\eps ^\delta \} \,.\end{equation*}
Recall that $\Prob_{5/8}$ is the law of a chordal $SLE _{8/3}$
from 0 to $\infty$ in $\H$, henceforth, we will simply call it $\Prob$. In
order to make the formulas more concise we will denote the event
$\{ z \mbox{ is to the right of $\gamma$ } \}$ by $R(z)$. Then,
\begin{equation}\label{kk} \Prob \{ R(\text{m}_\eps ^{-1}(\z)) |(F_
\eps)^c \cap F_\eps ^\delta \} = \frac{\Prob \{ R(\text{m}_\eps
^{-1}(\z)) |  F_ \eps ^\delta \} \Prob \{ F _\eps ^\delta \} -
\Prob \{ R(\text{m}_ \eps ^ {-1}(\z)) |  F_ \eps \} \Prob \{ F_ \eps \}}
{\Prob \{ F _ \eps ^\delta \} - \Prob \{ F_ \eps \}}\,.
\end{equation} The reason for this last step is that now all the
probabilities involved can be computed using the {\em restriction}
property for $SLE _{8/3}$, and a simple formula, see lemma
\ref{l.1}, for the probability that a point is to the right of an
$SLE_{8/3}$ path from 0 to $\infty$ in $\H$. This requires (cf.
section~\ref{preliminaries}) to know the unique conformal map
$\Phi_\eps=\Phi_{D_\eps}$ from $\H\setminus D_\eps$ into $\H$,
with $\Phi_\eps(0)=0$, $\Phi_\eps(\infty)=\infty$ and
$\Phi'_\eps(\infty)=1$ (with a similar statement for $D_\eps^\delta$).
Thus by restriction, the law of the chordal $SLE_{8/3}$ in $\H$
conditioned not to touch $D_\eps$ is the inverse image of the
chordal $SLE$ in $\H$ by $\Phi_\eps$. This implies for the
quantities we need to compute :
\begin{eqnarray*}
\Prob\{R(\mei(\z))|F_\eps\}& = &\Prob\{R(\mei(\Phi_\eps(\z)))\}\\
\Prob\{R(\mei(\z))|F_\eps^\delta\}& =&
\Prob\{R(\mei(\Phi_\eps^\delta(\z)))\}\,.
\end{eqnarray*}
Note that $\mei$ is a M\"obius transformation, which maps $\infty$ to $-1$.
Therefore, $D_\eps$ and $D_\eps^\delta$ are half disks whose centers are very
close to -1.
The fact that they are not exactly centered at -1 is due to the
lack of symmetry in the problem : an $SLE$ from 0 to $\eps$ in a half
disk $\Dp$ centered in 0. Nevertheless, for the computation of
$\Phi_\eps(z)$ and $\Phi_\eps^\delta(z)$, we can think of $D_\eps$
and $D_\eps^\delta$ as two half disks centered at -1 with radii
respectively $\eps$ and $(1-\delta)\eps$. If we carried out  the
computations with the actual disks (straightforward but tedious),
we would see that our
approximation is of order
$O(\eps^2 + \eps ^2 \delta^2/ |z + 1| + \eps^4/|z+1|^2)$, when $z$
goes to -1. In this way, we have
\begin{eqnarray*}
&\Phi_\eps(z)&=z-\eps^2+\frac{\eps^2}{z+1}+O(\eps^2+\f {\eps^4}{|z+1|^2})
\\ &\Phi_\eps^\delta(z)&=
z-\eps^2(1-\delta)^2+\frac{\eps^2(1-\delta)^2}{z+1}+
O(\eps^2+\f{\eps^2 \delta^2 }{ |z +1| }+ \f{\eps^4}{|z+1|^2})\,.
\end{eqnarray*}
We now have to evaluate these functions at the point
$\mei(\z)=\z/(\eps-\z)=-1-\frac{\eps}{\z}+O(\eps^2)$ (recall $\z$
is fixed). The approximations $O(\eps^4/|z+1|^2)$ and $O(\eps^2
\delta ^2 / |z + 1|)$ at the point $\mei(\z)$ are of order
$O(\eps^2)$ and $O(\eps \delta^2)$, respectively; this gives us :
\begin{eqnarray*}
\Phi_\eps^\delta(\mei(\z))&=& -1-\frac{\eps}{\z}+
\frac{\eps^2(1-\delta)^2}{-\eps/\z+O(\eps^2)}+O(\eps \delta ^2 + \eps^2)\\
&=&-1-\eps(\z+\frac{1}{\z})+2\eps\delta\z+O(\eps\delta^2+\eps^2)\,.
\end{eqnarray*}
Using the Taylor series for the logarithm, and then taking the
imaginary part, we see that
\begin{equation*}
\arg\left(  \Phi _\eps ^\delta (\mei(\z)) \right) = \pi + \eps
{\rm Im} (\z + \f 1 \z ) - 2 \eps \delta {\rm Im} (\z) +
O(\eps\delta ^2 + \eps ^2).
\end{equation*}
Now, using lemma \ref{l.1}, and the Taylor series for cosine we
see that
\begin{equation}\label{pro}
\Prob \{ R(\Phi _\eps ^\delta (\mei(\z)) ) \} = \f {\eps ^2} 4
\left[ \left({\rm Im} (\z + \f 1 \z)\right)^2 -4 \delta {\rm Im}
(\z + \f 1 \z ) {\rm Im} (\z) \right] + O(\eps^2\delta ^2 + \eps
^3).
\end{equation}
In particular, if we set $\delta = 0$ we obtain,
\begin{equation}\label{kkk}
\Prob \{ R(\Phi _\eps (\mei(\z)) ) \} = \f {\eps ^2} 4 \left({\rm
Im} (\z + \f 1 \z)\right)^2 + O(\eps ^3).
\end{equation}
Also, by \eqref{restriction}, we have (our approximation
doesn't change significantly the derivative at 0 which is far away
from small disks centered at -1) :
\begin{eqnarray*}
\Prob\{F_\eps^\delta\}=\Prob_{5/8}\{\gamma[0,\infty)\cap
D_\eps^\delta=\emptyset\}&=&(\Phi_\eps^\delta)'(0)^{5/8}\\
&=&(1-\eps^2(1-2\delta+O(\delta^2)))^{5/8}+O(\eps^3)\\
&=&1-\frac{5}{8}\eps^2+\frac{5}{4}\eps^2\delta+O(\eps^2\delta^2+\eps^3)
\end{eqnarray*}
Similarly, $\Prob\{F_\eps\}=1-5/8\eps^2+O(\eps^3)$, which gives
\begin{equation}
\label{remark}
\Prob\{F_\eps^\delta\}-\Prob\{F_\eps\}=\frac{5}{4}\eps^2\delta+O(\eps^2\delta^2+\eps^3)\,.
\end{equation}

Hence, by combining this last expression, \ref{kk}, \ref{pro},
\ref{kkk} and using the fact that both $\Prob\{F_\eps\}$ and
$\Prob\{F_\eps^\delta\}$ are $1+O(\eps^2)$, we obtain :
\begin{align*}
\lim _{\delta \downarrow 0} \lim _{\eps \downarrow 0} \Prob_\eps
\{ \z \mbox{ inside } & | \{\gs \in [1, 1 + \delta) \}\} \\
 &= 
\lim _{\delta \downarrow 0} \lim _{\eps \downarrow 0} \f {\f
{\eps ^2} 4 (-4 \delta  {\rm Im} (\z + \f 1 \z ) {\rm Im} (\z) ) +
\eps ^2 O(\delta ^2) + O(\eps ^3) }
{\f 5 {4} \eps ^2 \delta + O(\eps^ 2 \delta ^2 + \eps ^3)}\\
&= - \f 4 5 {\rm Im}(\z + \f 1 \z ) {\rm Im} (\z).
\end{align*}

Therefore, by \eqref{goal}, and using polar coordinates to
evaluate the integral, we get :
\begin{eqnarray*}
\slemeas( A |\gs = 1) &=& \int_{\Dp}
- \f 4 5 {\rm Im}(z + \f 1 z ) {\rm Im} (z) d{\rm A}(z)\\
&=& \f \pi {10} \,.\end{eqnarray*}
This concludes the proof of the lemma.
\end{proof}

Below we state a result that we have used extensively in our proof; for the reader's sake
we will sketch a proof.
This lemma gives an equivalent expression to the one given by
Schramm, in \cite{schramm}.
\begin{lemma}\label{l.1} Let $\gamma$ be chordal $SLE _{\kappa}$ in $\H$ with $\kappa\leq 4$, and
let $z = re ^{i \theta}$ be a point in $\H$. If we let $f(z) =
\Prob \{ z \mbox{ is to the right of } \gamma[0,\infty) \},$ then
$f$, which by scaling depends only on $\theta$, is given by
\[
f(\theta) = \frac{1}{\int _{0} ^\pi (\sin u)^{\f {2(4 -
\kappa)}{\kappa}} du}\;\int _{\theta} ^\pi (\sin u)^{\f {2(4 -
\kappa)}{\kappa}} du .\] In particular, for $\kappa = \f 8 3$ :
\[
\Prob \{ z \mbox{ is to the right of } \gamma[0,\infty) \} = 1/2 +
1/2 \cos (\theta).
\]
\end{lemma}
\begin{proof} {\em (sketch)}

As already mentioned, by scale invariance of $SLE$, the probability
that a point $z=r e^{i\theta}$ is to the right of the curve only depends
on the angle $\theta$. Thus, this probability is a certain function $f$ of the angle $\theta$. $SLE$
curves satisfy also a conformal-type Markov property. Thus, if $X_t$ is the
unique conformal map from $\H\setminus\gamma(0,t]$ onto $\H$ satisfying $X_t(\infty)=\infty$
, $X_t'(\infty)=1$ and $X_t(\gamma(t))=0$, we get :
\begin{equation*}
\Prob(z \mbox{ is on the right }|\mathcal{F}_t) = \Prob(X_t(z) \mbox{ is on the right}) =f(\theta_t)\,,
\end{equation*}
where $\theta_t$ is the continuous argument of $X_t(z)$. This shows that $f(\theta_t)_{t\geq 0}$ is a
martingale in $(0,1)$. Using that
$X_t(z)=g_t(z)-\sqrt{\kappa}dB_t$, where $g _t$ is defined by :
\begin{equation*}
\partial_t g_t(z)=\frac{2}{g_t(z)-\sqrt{\kappa}B_t}\,, \quad g_0(z)=z \qquad z\text{ in }\H,
\end{equation*}
we have
\begin{eqnarray*}
dX_t& =& dg_t(z)-\sqrt{\kappa}dB_t = \frac{2}{X_t}dt -\sqrt{\kappa}dB_t\,,\\
d\log X_t&=&\frac{2}{X_t^2}dt-\frac{\sqrt{\kappa}}{X_t}dB_t-\frac{\kappa}{2X_t^2}dt =\frac{(4-\kappa)}{2X_t^2}
dt-\frac{\sqrt{\kappa}}{X_t}dB_t\,,
\end{eqnarray*}
and by taking the imaginary part :
\begin{equation*}
d\theta_t=\frac{\kappa-4}{2|X_t|^2}\sin(2\theta_t)dt+\frac{\sqrt{\kappa}}{|X_t|}\sin(\theta_t)dB_t\,.
\end{equation*}
Now suppose $f$ is a $\mathcal{C}^2$ function, and apply It\^{o}'s formula to $f(\theta_t)$. We want this
process to be a martingale, so the $dt$ term in the expression for $df(\theta_t)$ has to be 0.
This gives a simple second order deterministic differential equation. Moreover we have the boundary conditions
$f(0)=1$ and $f(\pi)=0$. There is a unique solution, given in the lemma, which indeed is $\mathcal{C}^2$.
\end{proof}

{\bf Remark :} We would like to point out that the 1/5 in the
final result, comes from the 8/5 in the {\em restriction} formula
\ref{restriction}.

\section{Decomposition of the expected area of the Brownian loop into the expected
areas of the regions with fixed winding number}

Let $z\in\C\setminus \{0\}$ be fixed, 
and $(B_t)_{0\leq t\leq 1}$ a Brownian loop in $\C$ starting at
0. Almost surely $z \not \in \{B_s: 0 \le s \le 1  \}$, and therefore
we can define its index $n_z$.  
More precisely, $\forall s\in[0,1], B_s-z=R_s^z\exp(i \theta_s^z)$, where
$R_s^z=|B_s-z|$ and $\theta_s^z$ is any continuous 
representative of the argument. The index $n_z$ is by definition
$\frac{\theta_1^z-\theta_0^z}{2\pi}$ ; this is the number of times that 
the Brownian particle winds around $z$. For each $n\in\Z$, $n\neq 0$, 
let $\mathcal{W}_n$ denote the area
of the open set of points of index $n_z=n$. This random
variable can be written as :
\begin{equation*}
\mathcal{W}_n=\int_{\C} 1_{\{n_z=n\}}d{\rm A}(z)\,.
\end{equation*}
Let $\mathcal{W}_0$ be the area of the open set of points inside
the loop that have index zero :
\begin{equation*}
\mathcal{W}_0=\int_{\C} 1_{\{n_z=0\}\cap\{z\text{ is inside}\}}d{\rm A}(z)\,.
\end{equation*}
\begin{figure}
\begin{center}
\includegraphics[width=0.65\textwidth]{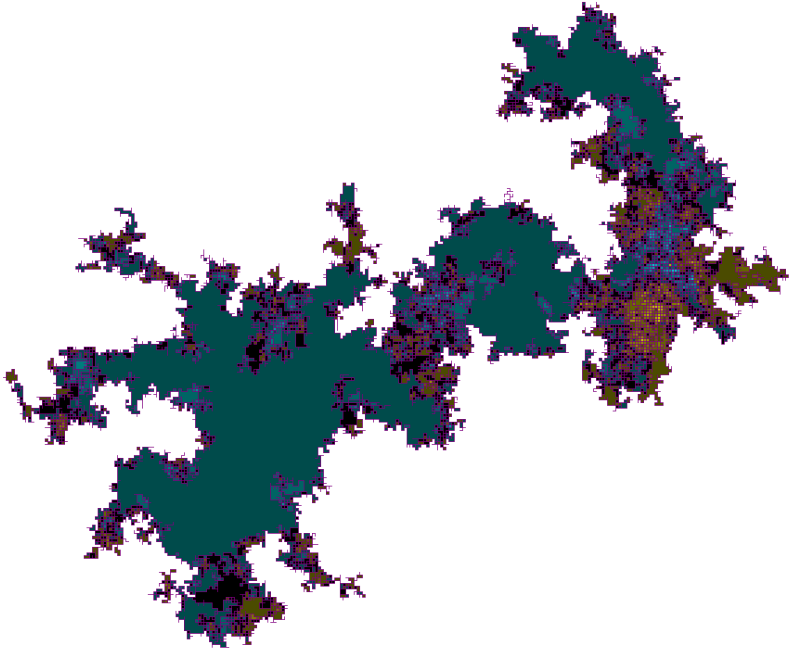}
\end{center}
\caption{Different indices in a random walk of 50000 steps, 
black areas correspond to index~0.}
\end{figure}
Since the Brownian curve is of Lebesgue measure zero, 
we have the following decomposition of the area $\A$ 
inside the Brownian loop (basically, the Brownian path
does not take much place inside its hull)
\begin{equation*}
\A=\sum_{n\in\Z} \mathcal{W}_n
\end{equation*}
Hence :
\begin{equation*}
\label{decomp.1}
\E(\A)=\frac{\pi}{5}=\sum_{n\in\Z}\E(\mathcal{W}_n).
\end{equation*}
Using results of Yor \cite{yor}, it will be straightforward to compute
$\E(\mathcal{W}_n)$ for $n\neq0$. And, hence, by subtracting 
from $\pi / 5$, one can obtain the value of $\E(\mathcal{W}_0)$.

\begin{theorem}

\begin{equation}
\E(\mathcal{W}_n) = \left\{ \begin{array}{ll} \frac{\pi}{30} \qquad &  n=0 \,, \\[3pt]
                       \frac{1}{2\pi n^2} \qquad & n\neq
                       0,\,n\in\Z\,.
                       \end{array}
                       \right.
                       \end{equation}

\end{theorem}

{\bf Remark :} This result is consistent with the asymptotic result
obtained by Werner in \cite{tournebeaucoup}, about the area
$A_n^t$ of the set of points around which the planar Brownian
motion (not the loop) winds around $n$ times on $[0,t]$. It is
indeed proved that $A_n^t$ is equivalent (in the $L^2$-sense) to
$\frac{t}{2\pi n^2}$ as $n$ goes to infinity. Very roughly 
the area of the $n$-sector for large $n$ 
comes from local contributions along the path, hence the 
global picture of the hull is not relevant; that is why, 
both Brownian motion and Brownian bridge should have 
the same asymptotics. Werner's proof requires to
compute the asymptotics of the first and second moments. This
present paper gives exact computations for the first moments in the
case of the loop, but it does not provide 
any information about the second moments.
\begin{proof}
We start by computing $\E(\mathcal{W}_n)$ for
$n\neq 0$. For this purpose we use theorem \ref{loi}, which was 
proved by Yor \cite{yor}.
Thus, for each $n\neq0$, using polar coordinates :
\begin{eqnarray*}
\E(\mathcal{W}_n) &=& \int_\C \Prob(n_z=n) d{\rm A}(z) \\
&=& 2 \pi \int_0^\infty r dr 
e^{-r^2}\left[\int_0^\infty dt
e^{-r^2\cosh(t)}\left(
\frac{2n - 1}{t^2+(2n-1)^2\pi^2}-\frac{2n + 1}{t^2+(2n+1)^2\pi^2}\right)\right]\\
&=&2\pi\int_0^\infty
dt\left(\frac{2n - 1}{t^2+(2n-1)^2\pi^2}-\frac{2n + 1}{t^2+(2n+1)^2\pi^2}\right)\int_0^\infty
r e^{-r^2(1+\cosh(t))}dr\\
&=&\pi\int_0^\infty
\frac{dt}{1+\cosh(t)}\left(\frac{2n - 1}{t^2+(2n-1)^2\pi^2}-
\frac{2n + 1}{t^2+(2n+1)^2\pi^2}\right) \\
&=& \frac 1 {2 \pi n^2}\, .
\end{eqnarray*}

We sketch one possible way to see how to obtain the last line
in the above chain of equalities. 

It is slightly more convenient to generalize a bit, so thinking of 
$2n$ as $x$ and using the symmetry of the integrand, we 
consider the function 
\begin{equation*}
F(x) = \int_{-\infty}^\infty
\frac{dt}{1+\cosh(t)}\left(\frac{x - 1}{t^2+(x-1)^2\pi^2}-
\frac{x + 1}{t^2+(x+1)^2\pi^2}\right).
\end{equation*}
In this new notation what we want to prove is that 
$F(x) = \frac 4 {\pi ^2 x ^2}$ (for $x \ge |2|$). Since, 
$F$ is symmetric about 0, it is enough to study the 
case of $x$ positive; furthermore, since $F$ is 
real analytic on $\{ x : x > 1 \}$, we can allow 
ourselves to assume that $x$ is not an integer.
Now, for $x > 1$ and $x$ not an integer,
a simple residue computation with appropriate contours yields
\begin{equation*}
F(x) = -\frac 8 {\pi ^2} \sum\limits_{k = 1}^\infty
\left(\frac{(2k - 1)(x-1)}{((x-1)^2 - (2k-1)^2)^2} - 
\frac{(2k-1)(x+1)}{((x+1)^2 - (2k - 1)^2)^2}\right).
\end{equation*}
In order to evaluate this sum, it is enough to 
notice that using partial fractions one can obtain
\begin{equation*}
\sum\limits_{k=1}^\infty
\f{(2k-1)w}{(w^2 -(2k-1)^2)^2} =
-\frac 1 {16} \sum\limits_{k=0}^\infty
\left( \f 1 {(k + w/2 + 1/2)^2} - 
\f 1 {(k - w/2 + 1/2)^2}\right),
\end{equation*}
and substituting $x -1$ and $x + 1$ for $w$, and noticing the 
telescoping cancellations one readily obtains 
\begin{equation*}
F(x) = \f 4 {\pi ^2 x ^2},
\end{equation*}
hence, 
$\E(\mathcal{W}_n)=\frac{1}{2\pi n^2}$.

Finally, using the fact that $\sum_{n=1}^\infty
\frac{1}{n^2}=\frac{\pi^2}{6}$, and the fact that 
the area of the Brownian loop is $\pi / 5$ we conclude
$\E(\mathcal{W}_0)=\frac {\pi} {30}$. This finishes the 
proof of the theorem.
\end{proof}

\begin{theorem}\label{loi}
Fix $z=r e^{i\theta}$, with $r \ne 0$. 
Under the law of a Brownian loop of time
duration one, starting at 0, we have the following probabilities :
\begin{eqnarray}
\Prob(n_z=n)&=&e^{-r^2}[\Psi_r((2n-1)\pi)-\Psi_r((2n+1)\pi)]\text{    if }n\in\Z\setminus 0\,, \\
\Prob(n_z=0)&=&1+e^{-r^2}[\Psi_r(-\pi)-\Psi_r(\pi)]\,,
\end{eqnarray}
where $\forall x\neq  0$,
\begin{equation*}
\Psi_r(x)=\frac{x}{\pi}\int_0^{\infty}e^{-r^2\cosh(t)}
\frac{dt}{t^2+x^2}\,.
\end{equation*}
\end{theorem}

 {\bf Acknowledgments :} We wish to thank Greg
Lawler and Wendelin Werner for suggesting the problem and for
fruitful discussions, and Wendelin Werner for pointing out the
link with the paper of Yor \cite{yor}.

\end{document}